# Comments on "Structure design of two types of sliding-mode controllers for a class of under-actuated mechanical systems"

B.L.Ma

The Seventh Research Division, Beijing University of Aeronautics and Astronautics, Beijing 100083

Email: ma_baoli@yahoo.com

**Abstracts:** In this note, it is shown by counterexamples, theoretical analysis and simulation tests that the two types of sliding-mode controllers presented in the paper [1] fail to solve the stabilization problem of a class of under-actuated mechanical systems.

## 1   Introduction

In the paper [1], two types of sliding-mode control approaches are proposed for the following under-actuated system

$$\begin{aligned}
\dot{x}_1 &= x_2 \\
\dot{x}_2 &= f_1(X) + b_1(X)u \\
\dot{x}_3 &= x_4 \\
\dot{x}_4 &= f_2(X) + b_2(X)u \\
&\vdots \\
\dot{x}_{2n-1} &= x_{2n} \\
\dot{x}_{2n} &= f_n(X) + b_n(X)u
\end{aligned} \qquad (a)$$

Under the simple hypothesis that:

(A1).  $0 \leq |f_i(X)| \leq M_i, X \in A_d^c$,

(A2).  $0 < |b_i(X)| \leq B_i, X \in A_d^c$,

the authors claim to achieve convergence of state errors to zero using the proposed two types of sliding mode controllers.

The results are a generalization of the approaches presented in [2]. We had pointed out that the approaches proposed in [2] are false in the previous comment [3]. In this note, it is indicated that the results presented in [1] are also false. We show this first via intuitive counterexamples and then indicate the inconsistencies in the proof given in [1] of Theorem 1 and Theorem 2, the main results of that paper. Simulation results and detail theoretical analysis for the model of overhead crane are also given to show the failures and useless of the proposed schemes presented in that paper.

## 2. Counterexamples for Theorem 1 and Theorem 2

Theorem 1 in [1] claims that: the proposed **IHSSMC** scheme guarantee that the closed-loop system is globally stable in the sense that all signals involved are bounded, with the errors converging to zero asymptotically.

Theorem 2 in [1] claims that: the proposed **AHSSMC** scheme guarantee that the closed-loop systems is globally stable in the sense that all signals involved are bounded, with the errors converging to zero asymptotically.

The following counterexamples indicate that the above two claims are both false.

Consider the four states case of system (a)

$$\begin{aligned}\dot{x}_1 &= x_2 \\ \dot{x}_2 &= f_1 + b_1 u \\ \dot{x}_3 &= x_4 \\ \dot{x}_4 &= f_2 + b_2 u\end{aligned} \quad (b)$$

Assume that $f_1 = kf_2, b_1 = kb_2 \neq 0$ ($k$ is a nonzero constant), it is obvious that system (b) satisfy the hypothesis (A1) and (A2), however, system (b) can not be stabilized to origin by any control input for initial states satisfying $x_1(0) \neq 0, x_2(0) = x_3(0) = x_4(0) = 0$. The analysis is as follows: $f_1 = kf_2, b_1 = kb_2 \neq 0$ imply $\dot{x}_2(t) - k\dot{x}_4(t) = 0 (\forall t \geq 0)$, which in tune imply $x_2(t) - kx_4(t) = x_2(0) - kx_4(0) = 0 (\forall t \geq 0)$, so that $x_1(t) - kx_3(t) = x_1(0) - kx_3(0) = x_1(0) \neq 0 (\forall t \geq 0)$, i.e. $x_1(t) - kx_3(t)$ keeps to be a nonzero constant $x_1(0)$ for all $t \geq 0$, thus $x_1(t), x_3(t)$ can not be simultaneously zero for all $t \geq 0$.

The above example is not pure artificial, it can be found in real robot systems. To show this, let us consider a Pendubot shown in Fig.a, a well-known under-actuated planar manipulator with the first joint actuated and the second one free.

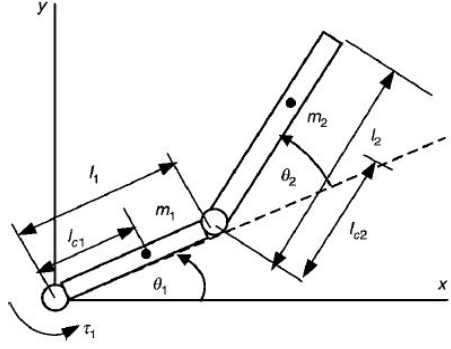

Fig.a. The Pendubot:two-link manipulator with a free second joint.

The dynamic model of Pendubot can be obtained as [2]:

$$D(\theta)\ddot{\theta} + C(\theta,\dot{\theta})\dot{\theta} + G(\theta) = \tau$$

where

$$\theta = \begin{bmatrix} \theta_1 \\ \theta_2 \end{bmatrix}, \tau = \begin{bmatrix} \tau_1 \\ 0 \end{bmatrix},$$

$$D(\theta) = \begin{bmatrix} q_1 + q_2 + 2q_3 \cos\theta_2 & q_2 + q_3 \cos\theta_2 \\ q_2 + q_3 \cos\theta_2 & q_2 \end{bmatrix},$$

$$C(\theta,\dot{\theta}) = q_3 \sin(\theta_2) \begin{bmatrix} -\dot{\theta}_2 & -\dot{\theta}_2 - \dot{\theta}_1 \\ \dot{\theta}_1 & 0 \end{bmatrix},$$

$$G(\theta) = \begin{bmatrix} q_4 g \cos\theta_1 + q_5 g \cos(\theta_1 + \theta_2) \\ q_5 g \cos(\theta_1 + \theta_2) \end{bmatrix},$$

$$q_1 = m_1 l_{c1}^2 + m_2 l_1^2 + I_1, q_2 = m_2 l_{c2}^2 + I_2,$$

$$q_3 = m_2 l_1 l_{c2}, q_4 = m_1 l_{c1} + m_2 l_1, q_5 = m_2 l_{c2}$$

Now consider the special case with $l_{c2} = 0$, i.e. the center of mass of second-link lies on the second-joint axis. In this case, we have $q_3 = 0, q_5 = 0$, thus

$$D(\theta) = \begin{bmatrix} q_1 + q_2 & q_2 \\ q_2 & q_2 \end{bmatrix}, C(\theta,\dot{\theta}) = 0, G(\theta) = \begin{bmatrix} q_4 g \cos\theta_1 \\ 0 \end{bmatrix}$$

The state space model in this special case can be written as

$$\begin{bmatrix} \ddot{\theta}_1 \\ \ddot{\theta}_2 \end{bmatrix} = D^{-1}[\tau - C(\theta,\dot{\theta})\dot{\theta} - G(\theta)] = \frac{1}{q_1}\begin{bmatrix} \tau_1 - q_4 g \cos\theta_1 \\ -(\tau_1 - q_4 g \cos\theta_1) \end{bmatrix} = \begin{bmatrix} u \\ -u \end{bmatrix} \quad (c)$$

where $u \triangleq \frac{1}{q_1}(\tau_1 - q_4 g \cos\theta_1)$.

It is clear that system (c) is just in the form of (b) with $f_1 = f_2 = 0, b_1 = -b_2 = 1 \neq 0$.

## 3. Errors in the Proof

Theorem 1 and Theorem 2 in [1] provide the theoretical basis for the authors' claims that the proposed sliding mode controllers achieve state regulation. We have shown that the claims of Theorem 1 and Theorem 2 are both false by intuitive counterexamples. In what follows, some serious inconsistencies in the proof of Theorem 1 and Theorem 2 will be indicated.

First of all, note that the proof of Theorem 1 relies on Eq.(31)-Eq.(32), which are inferred from Eq.(30). But this is incorrect, since Eq.(30) only implies $\int_0^\infty c_i^2 x_{i+1}^2 dt < \infty$, $\int_0^\infty s_{i-1}^2 dt < \infty$, and $\int_0^\infty x_{i+1}^2 dt < \infty$ can not be deduced from $\int_0^\infty c_i^2 x_{i+1}^2 dt < \infty$ as $c_i = C_i sign(x_{i+1} s_{i-1}) = C_i sign(s_1) sign(x_{i+1})$ is zero for $s_1 = 0$. Similarly, only $\int_0^\infty |c_i x_{i+1}| dt < \infty$, $\int_0^\infty |s_{i-1}| dt < \infty$ can be deduced from Eq.(33). Hence, only $\lim_{t\to\infty} s_{i-1} = 0, \lim_{t\to\infty} c_i x_{i+1} = 0$ can be concluded by Barbalat lemma. As $\lim_{t\to\infty} s_1 = 0$, so $\lim_{t\to\infty} c_i = \lim_{t\to\infty} C_i sign(x_{i+1} s_{i-1}) = \lim_{t\to\infty} C_i sign(s_1) sign(x_{i+1}) = 0$, thus $\lim_{t\to\infty} x_{i+1} = 0$ can not be deduced from $\lim_{t\to\infty} c_i x_{i+1} = 0$.

In fact, according to the definition of $s_i$ and $c_i$, $s_{2n-1}$ can be expressed as

$$s_{2n-1} = s_1 + C_2 sign(s_1 x_3) x_3 + \cdots + C_{2n-1} sign(s_1 x_{2n}) x_{2n} = sign(s_1)\left(|s_1| + \sum_{i=3}^{2n} C_{i-1} |x_i|\right)$$, which

means that $s_{2n-1} = 0 \Leftrightarrow s_1 = 0$, i.e. the sliding surface $s_{2n-1} = 0$ is completely equivalent to the sliding surface $s_1 = 0$, therefore, the **IHSSMC** scheme can only guarantee that the state trajectories reach to the surface $s_1 = x_2 + c_1 x_1 = \dot{x}_1 + c_1 x_1 = 0$. In the sliding surface $s_1 = 0$, if $c_1$ is chosen as a positive constant, $\lim_{t\to\infty} x_1 = 0, \lim_{t\to\infty} x_2 = 0$ can be concluded, if $c_1$ is chosen as $c_1 = C_1 sign(x_1 x_2)$, only $\lim_{t\to\infty} x_2 = 0$ can be concluded since

$$s_1 = x_2 + c_1 x_1 = sign(x_2)(|x_2| + C_1 |x_1|) = 0 \Leftrightarrow x_2 = 0.$$

The proof of Theorem 2 in [1] relies on $s_i \in L_\infty$, $\dot{s}_i \in L_\infty$ and $s_i \in L_2$, which in tune rely on (59)-(61) and (65)-(66), but this is not true. First, $s_i \in L_\infty$, $\dot{s}_i \in L_\infty$ can not be deduced from (59)-(61) unless $\|X(t)\|_{\infty,w} \in L_\infty$ (or equivalently $X(t) \in A_d^c$) is guaranteed, but $X(t) \in A_d^c$ can not be guaranteed by **AHSSMC** scheme as the dynamics on the sliding surface $S = 0$ may be unstable with some or all states diverging to infinity. Secondly, $s_i \in L_2$ is deduced from (65)-(66), i.e. $S_1 \in L_2, S_2 \in L_2$, which is in tune deduced from $S \in L_2$, but this is not true since

$$S_1 = S + (\alpha_{i1} - \alpha_i)s_i \in L_2 \quad , \quad S_2 = S + (\alpha_{i2} - \alpha_i)s_i \in L_2 \quad \text{can not be inferred}$$

from $S \in L_2$ for $\alpha_{i1} \neq \alpha_{i2}$.

## 4. Verification and explanations of simulation results

The controllers presented in [1] had been used to stabilize an overhead crane. The simulation results in that paper show that the proposed sliding controller seems to work well: how can this take place?

### 4.1 Verification and explanations of simulation results with IHSSMC scheme

We feel that the simulation results of **IHSSMC** scheme presented in [1] are doubtful. To verify this, we have simulated the model of overhead crane using **IHSSMC** scheme with the same parameters as in [1]: $M = 1.0kg$, $m = 0.8kg$, $L = 0.305m$, $c_1 = C_1 = 1.4$, $C_2 = 0.2$, $C_3 = 0.1$, $k = 0.1$, $\eta = 1$. Our simulation results are shown in Fig.b. It is clear that the simulation results presented in [1] (Fig.2-3 in that paper) can not be completely reproduced. The time plots of $(x(t), \dot{x}(t))$ are similar, but the time plots of $(\theta(t), \dot{\theta}(t))$ are quite different. The simulation results presented in [1] shows that $(\theta(t), \dot{\theta}(t))$ converge to zero, ours show that $(\theta(t), \dot{\theta}(t))$ converge to a periodic trajectories. Our results are inconsistent with what the Theorem 1 in [1] predicted, but is consistent with our analysis, i.e. the **IHSSMC** scheme can only guarantee $\lim_{t \to \infty} s_1 = 0$ in finite time (as $c_1$ is chosen as positive constant, so $\lim_{t \to \infty} x_1 = 0, \lim_{t \to \infty} x_2 = 0$ are guaranteed).

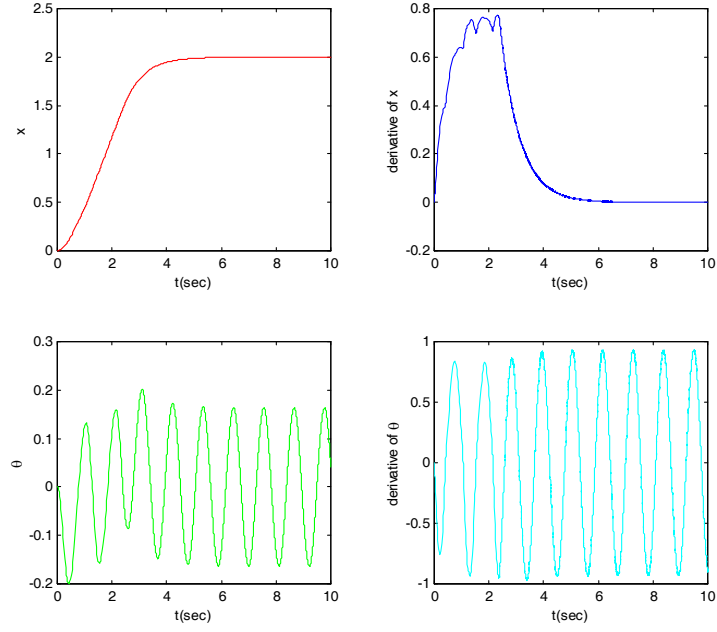

Fig.b Simulation results using **IHSSMC** scheme proposed in [1].

For the conveniences of reviewers to check and run the simulation programs, the simulation codes are copied in Appendix A.

In what follows, the periodic behavior of $(\theta, \dot{\theta})$ will be further analyzed and explained by exploration of dynamics on the sliding surface $s_1 = 0$.

Denote $X = col[x_1, x_2, x_3, x_3]$, $x_1 = x - x_d$, $x_2 = \dot{x}$, $x_3 = \theta$, $x_4 = \dot{\theta}$, the dynamics of overhead crane can be written as [1]:

$$\dot{x}_1 = x_2, \dot{x}_2 = f_1 + b_1 u, \dot{x}_3 = x_4, \dot{x}_4 = f_2 + b_2 u \tag{d}$$

where

$$f_1 = \frac{mLx_4^2 \sin x_3 + mg \sin x_3 \cos x_3}{M + m\sin^2 x_3}, \qquad b_1 = \frac{1}{M + m\sin^2 x_3},$$

$$f_2 = -\frac{(m+M)g \sin x_3 + mLx_4^2 \sin x_3 \cos x_3}{(M + m\sin^2 x_3)L}, \qquad b_2 = -\frac{\cos x_3}{(M + m\sin^2 x_3)L} \tag{e}$$

As discussed in Section 3, the **IHSSMC** scheme can only guarantee the state trajectories reach to the surface $s_1 = x_2 + c_1 x_1 = 0$ in finite time, and then $(x_1, x_2)$ will converge to zero along with $s_1 = 0$.

The crane dynamics on the sliding surface $s_1 = 0$ can be analyzed as follows.

By $s_1 = x_2 + c_1 x_1 = 0$ and $\dot{s}_1 = c_1 x_2 + f_1 + b_1 u = 0$, one gets $x_2 = -cx_1$ and $u = u_{eq} \triangleq \dfrac{c_1^2 x_1 - f_1(x_1, -cx_1, x_3, x_4)}{b_1(x_1, -cx_1, x_3, x_4)}$. Substituting $x_2 = -c_1 x_1, u = u_{eq}$ into the dynamics of $x_3, x_4$ results

$$\dot{x}_3 = x_4,$$
$$\dot{x}_4 = f_2(x_1, -cx_1, x_3, x_4) + b_2(x_1, -cx_1, x_3, x_4) u_{eq}$$
$$= f_2(x_1, -cx_1, x_3, x_4) + \frac{b_2(x_1, -cx_1, x_3, x_4)}{b_1(x_1, -cx_1, x_3, x_4)} \left(c_1^2 x_1 - f_1(x_1, -cx_1, x_3, x_4)\right)$$
$$= -\frac{(m+M)g \sin x_3 + mLx_4^2 \sin x_3 \cos x_3}{(M + m \sin^2 x_3) L} - \frac{\cos x_3}{L}\left(c_1^2 x_1 - \frac{mLx_4^2 \sin x_3 + mg \sin x_3 \cos x_3}{M + m \sin^2 x_3}\right)$$
$$= -\frac{g}{L} \sin x_3 - \frac{c_1^2}{L} x_1 \cos x_3$$

As $\lim\limits_{t \to \infty} x_1 = 0$, so the asymptotical dynamics of $x_3, x_4$ becomes

$$\dot{x}_3 = x_4, \dot{x}_4 = -\frac{g}{L} \sin x_3$$

Take $V = 0.5 x_4^2 + g/L(1 - \cos x_3)$, then $\dot{V} = 0$, so that the asymptotical dynamics of $x_3, x_4$ is Lyapunov stable and the trajectories of $x_3, x_4$ converge to periodic ones!

To further show the invalidity and useless of **IHSSMC** scheme presented in [1], we design a linear control law to stabilize the overhead crane. The linearization of overhead crane dynamics at the desired equilibrium state $(x_d, 0, 0, 0)$ can be obtained as

$$\dot{X} = A_1 X + Bu$$

where $A_1 = \begin{bmatrix} 0 & 1 & 0 & 0 \\ 0 & 0 & \dfrac{mg}{M} & 0 \\ 0 & 0 & 0 & 1 \\ 0 & 0 & -\dfrac{(m+M)g}{ML} & 0 \end{bmatrix}, B = \begin{bmatrix} 0 \\ \dfrac{1}{M} \\ 0 \\ -\dfrac{1}{ML} \end{bmatrix}$.

It can be easily verified $\{A_1, B\}$ is controllable. Choose the desired poles as $P = (-3, -2.8, -2.6, -2.4)$, the feedback gain can be obtained as $K = acker(A_1, B, P)$ $= [1.3051\ 1.9468\ 7.3103\ -2.1602]$. Simulation results by applying linear control law $u = -K[x_1 - 2, x_2, x_3, x_4]^T$ are shown in Fig.c. It is clear that all the states converge to their desired values.

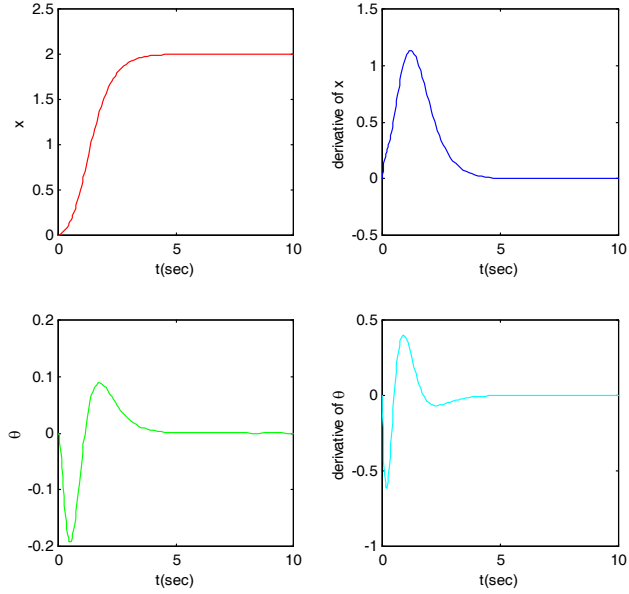

Fig.c Simulation results using linearization method.

**4.2 Verification and explanations of simulation results with AHSSMC scheme**

The simulation results of **AHSSMC** scheme presented in [1] are also doubtful. To verify this, we have simulated the overhead crane model using **AHSSMC** scheme proposed in that paper with the same parameters: $M = 1.0kg$, $m = 0.8kg$, $L = 0.305m$, $c_1 = 0.8$, $c_2 = 35$, $\alpha_1 = 10$, $\alpha_2 = 1$, $\eta = 3.5$, $k = 6$. Our simulation results are shown in Fig.d, indicating that all the state trajectories diverge to infinity (although the sliding variable $S$ converge to zero in finite time), while the simulation results in [1] (Fig.8-9 in that paper) shows that all the states converge to their desired values. Our results are inconsistent with what the Theorem 2 in [1] predicted, but are consistent with our analysis, i.e. the **AHSSMC** scheme can only guarantee finite convergence of $S$, but not all states, since the dynamics on the sliding surface $S = 0$ is not designed to be asymptotically stable.

For the convenience of reviewers to check and run the simulation programs, the simulation codes are copied in Appendix B.

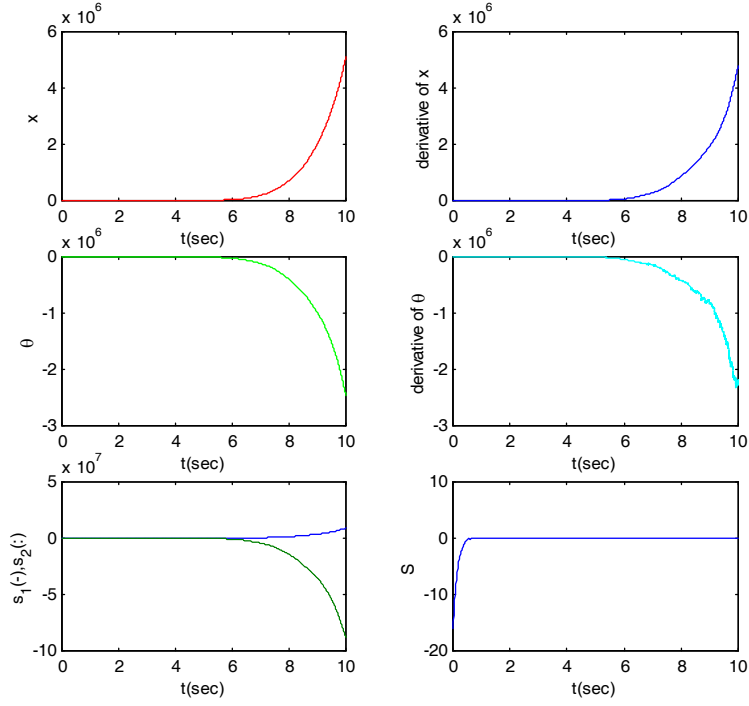

Fig.d. Simulation results using **AHSSMC** scheme proposed in [1].

The divergence of state trajectories can be further analyzed and explained by exploration of dynamics on sliding surface $S = 0$ as follows.

Set $\alpha_2 = 1$, the sliding surface of **AHSSMC** scheme can be expressed as

$$S = \alpha_1 s_1 + s_2 = \alpha_1(x_2 + c_1 x_1) + (x_4 + c_2 x_3) = \alpha_1 x_2 + \alpha_1 c_1 x_1 + x_4 + c_2 x_3 = 0$$

In the sliding surface $S = 0$, we have

$$x_4 = -c_2 x_3 - \alpha_1(x_2 + c_1 x_1) \tag{f}$$

The derivative of $S$ can be calculated as

$$\begin{aligned}\dot{S} &= \frac{d}{dt}(\alpha_1 x_2 + \alpha_1 c_1 x_1 + x_4 + c_2 x_3) \\ &= \alpha_1 f_1 + \alpha_1 b_1 u + c_1 \alpha_1 x_2 + f_2 + b_2 u + c_2 x_4) \\ &= \alpha_1 f_1 + \alpha_1 c_1 x_2 + f_2 + c_2 x_4 + (b_2 + \alpha_1 b_1)u \\ &= \alpha_1 f_1 + \alpha_1 c_1 x_2 + f_2 + c_2(-c_2 x_3 - \alpha_1(x_2 + c_1 x_1)) + (b_2 + \alpha_1 b_1)u \\ &= \alpha_1 f_1 + f_2 - c_1 c_2 \alpha_1 x_1 + \alpha_1(c_1 - c_2)x_2 - c_2^2 x_3 + (b_2 + \alpha_1 b_1)u\end{aligned}$$

From $\dot{S} = 0$, the equivalent control input can be solved as

$$u = u_{eq} \triangleq -\frac{1}{b_2 + \alpha_1 b_1}\left[\alpha_1 f_1 + f_2 - c_1 c_2 \alpha_1 x_1 + \alpha_1(c_1 - c_2)x_2 - c_2^2 x_3\right] \tag{g}$$

Therefore, the dynamics restricted on $S = \dot{S} = 0$ can be obtained as

$$\begin{aligned}
\dot{x}_1 &= x_2, \\
\dot{x}_2 &= f_1(x_1, x_2, x_3, -c_2 x_3 - \alpha_1(x_2 + c_1 x_1)) + b_1(x_1, x_2, x_3, -c_2 x_3 - \alpha_1(x_2 + c_1 x_1))u_{eq} \\
\dot{x}_3 &= -c_2 x_3 - \alpha_1(x_2 + c_1 x_1)
\end{aligned} \quad (h)$$

After substitution of expressions of $f_1, f_2, b_1, b_2, u_{eq}$ into (h), the linearization of system (h) can be obtained as

$$\begin{aligned}
\dot{x}_1 &= x_2, \\
\dot{x}_2 &= a_1 x_3 - \frac{b_{10}}{b_{20} + \alpha_1 b_{10}}\left(\alpha_1 a_1 x_3 + a_2 x_3 - c_1 c_2 \alpha_1 x_1 + \alpha_1(c_1 - c_2)x_2 - c_2^2 x_3\right) \\
&= \frac{b_{10} c_1 c_2 \alpha_1}{b_{20} + \alpha_1 b_{10}} x_1 - \frac{b_{10}\alpha_1(c_1 - c_2)}{b_{20} + \alpha_1 b_{10}} x_2 + \left[a_1 - \frac{b_{10}}{b_{20} + \alpha_1 b_{10}}(a_2 + \alpha_1 a_1 - c_2^2)\right]x_3 \\
\dot{x}_3 &= -c_2 x_3 - \alpha_1(x_2 + c_1 x_1)
\end{aligned} \quad (i)$$

where $a_1 = \dfrac{mg}{M}, a_2 = -\dfrac{(m+M)g}{ML}, b_{10} = \dfrac{1}{M}, b_{20} = -\dfrac{1}{ML}$.

System (i) can be rewritten compactly as

$$\begin{bmatrix} \dot{x}_1 \\ \dot{x}_2 \\ \dot{x}_3 \end{bmatrix} = A_2 \begin{bmatrix} x_1 \\ x_2 \\ x_3 \end{bmatrix} \quad (j)$$

where

$$A_2 = \begin{bmatrix} 0 & 1 & 0 \\ \dfrac{b_{10}}{b_{20} + \alpha_1 b_{10}} c_1 c_2 \alpha_1 & -\dfrac{b_{10}}{b_{20} + \alpha_1 b_{10}}\alpha_1(c_1 - c_2) & a_1 - \dfrac{b_{10}}{b_{20} + \alpha_1 b_{10}}(a_2 + \alpha_1 a_1 - c_2^2) \\ -\alpha_1 c_1 & -\alpha_1 & -c_2 \end{bmatrix}$$

For the parameters chosen in [1]: $c_1 = 0.8, c_2 = 35, \alpha_1 = 10$, the eigenvalues of $A_2$ can be computed as $\lambda(A_2) = \{-0.6527, 11.3918, 5.1438\}$, indicating that there are two unstable eigenvalues, therefore the dynamics on the sliding surface is unstable. This explains why the state trajectories diverge!

To guarantee asymptotically convergence of state trajectories, the correct method to design sliding surface should be: select parameters $c_1, c_2, \alpha_1$ such that the crane dynamics restricted on $S = 0$ is asymptotically stable, i.e. $c_1, c_2, \alpha_1$ should be selected such that $A_2$ is Hurwitz.

Direct computation results

$$\det(sI - A_2) = s^3 + l_1 s^2 + l_2 s + l_3$$

where

$$l_1 = c_2 + \frac{b_{10}}{b_{20} + \alpha_1 b_{10}} \alpha_1 (c_1 - c_2),$$

$$l_2 = \frac{\alpha_1}{b_{20} + \alpha_1 b_{10}} (b_{20} a_1 - b_{10} a_2),$$

$$l_3 = \frac{c_1 \alpha_1}{b_{20} + \alpha_1 b_{10}} (b_{20} a_1 - b_{10} a_2) = c_1 l_2$$

For the given desired Hurwitz polynomial: $d(s) = s^3 + d_1 s^2 + d_2 s + d_3$, $c_1, c_2, \alpha_1$ should be selected such that

$$d_1 = c_2 + \frac{b_{10}}{b_{20} + \alpha_1 b_{10}} \alpha_1 (c_1 - c_2),$$

$$d_2 = \frac{\alpha_1}{b_{20} + \alpha_1 b_{10}} (b_{20} a_1 - b_{10} a_2), \qquad (k)$$

$$d_3 = \frac{c_1 \alpha_1}{b_{20} + \alpha_1 b_{10}} (b_{20} a_1 - b_{10} a_2) = c_1 d_2$$

The design parameters $(c_1, c_2, \alpha_1)$ can be solved from (k) as

$$\alpha_1 = \frac{b_{20} d_2}{(b_{20} a_1 - b_{10} a_2 - d_1 b_{10})}, \; c_1 = \frac{d_3}{d_2}, \; c_2 = \frac{d_1 (b_{20} + \alpha_1 b_{10}) - b_{10} \alpha_1 c_1}{b_{20}} = d_1 + \frac{\alpha_1 b_{10} (d_1 - c_1)}{b_{20}} \qquad (l)$$

For the desired eigenvalues $\lambda_{d1} = -5, \lambda_{d1} = -4, \lambda_{d1} = -3$, $d_1, d_2, d_3$ are obtained as $d_1 = 12, d_2 = 47, d_3 = 60$, and $c_1, c_2, \alpha_1$ are computed as $c_1 = 1.2766$, $c_2 = -21.8964$, $\alpha_1 = 10.3638$. The simulation results corresponding to these new design parameters are shown in Fig.e. It is seen that all the states converge to their desired values.

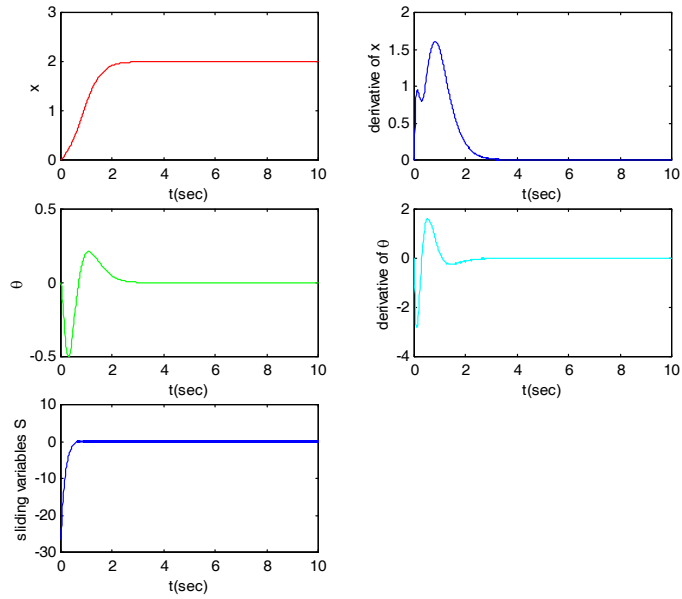

Fig.e Simulation results with stable sliding dynamics.

## 5. Conclusion

It has been shown by intuitive counterexamples, theoretical analysis and simulation tests that the sliding controllers presented in [1] are not only unable to solve the asymptotical stabilization problem of a class of under-actuated systems, but also unable to solve the asymptotical stabilization problem of simple linear (or linearizable) controllable systems.

**Acknowledgments**

This work was supported by National Science and Technology Major Project(No.2012CB821202), Beijing Natural Science Foundation(No.4122043), and National Nature Science Foundation (No.60874012).

**Appendix A. Simulation codes for stabilization of overhead crane using IHSSMC scheme proposed in [1].**

**%Main Program**

```
clear
global C1 C2 C3 eta k m M L g
C1=1.4; C2=0.2; C3=0.1; eta=1;k=0.1; m=0.8; M=1.0; L=0.305; g=9.8;
T=10; ts=[0,T];y0=[0;0;0;0];
tol1=1e-3; er=1e-4;
options =odeset('RelTol',tol1,'AbsTol',[er er er er]);
[t,y]=ode45('tem8',[0,T],y0,options);
subplot(2,2,1),plot(t,y(:,1),'r');xlabel('t(sec)');ylabel('x');
subplot(2,2,2),plot(t,y(:,2),'b');xlabel('t(sec)');ylabel('derivative of x');
subplot(2,2,3),plot(t,y(:,3),'g');xlabel('t(sec)');ylabel('\theta');
subplot(2,2,4),plot(t,y(:,4),'c');xlabel('t(sec)');ylabel('derivative of \theta');
```

**%Sub-Program**

```
function dy=tem8(t,y)
global C1 C2 C3 eta k m M L g
y1=y(1)-2;
c1=C1;s1=y(2)+c1*y1;
c2=C2*sign(y(3)*s1);s2=s1+c2*y(3);
c3=C3*sign(y(4)*s2);s3=s2+c3*y(4);
f1=(m*L*y(4)^2*sin(y(3))+m*g*sin(y(3))*cos(y(3)))/(M+m*(sin(y(3)))^2);
b1=1/(M+m*(sin(y(3)))^2);
f2=-((m+M)*g*sin(y(3))+m*L*y(4)^2*sin(y(3))*cos(y(3)))/((M+m*(sin(y(3)))^2)*L);
b2=-cos(y(3))/((M+m*(sin(y(3)))^2)*L);
ueq=-(c3*f2+c2*y(4)+f1+c1*y(2))/(c3*b2+b1);
```

```
usw=-(eta*sign(s3)+k*s3)/(c3*b2+b1);
u=ueq+usw;
dy=[y(2);f1+b1*u;y(4);f2+b2*u];
```

**Appendix B. Simulation codes for stabilization of overhead crane using AHSSMC scheme proposed in [1].**

```
%Main Program
clear
global c1 c2 alpha1 alpha2 eta k m M L g
m=0.8;M=1.0;L=0.305;g=9.8;
c1=0.8;c2=35;alpha1=10;alpha2=1;eta=3.5;k=6;
T=10;ts=[0,T];y0=[0;0;0;0];
tol1=1e-3;er=1e-4;
options = odeset('RelTol',tol1,'AbsTol',[er er er er]);
[t,y]=ode45('tem14',[0,T],y0,options);
y1=y(:,1)-2;
s1=y(:,2)+c1*y1;s2=y(:,4)+c2*y(:,3);
s=alpha1*s1+alpha2*s2;
subplot(3,2,1),plot(t,y(:,1),'r');xlabel('t(sec)');ylabel('x');
subplot(3,2,2),plot(t,y(:,2),'b');xlabel('t(sec)');ylabel('derivative of x');
subplot(3,2,3),plot(t,y(:,3),'g');xlabel('t(sec)');ylabel('\theta');
subplot(3,2,4),plot(t,y(:,4),'c');xlabel('t(sec)');ylabel('derivative of \theta');
subplot(3,2,5),plot(t,s1,'-',t,s2,':');xlabel('t(sec)');ylabel('s_1(-),s_2(:)');
subplot(3,2,6),plot(t,s);xlabel('t(sec)');ylabel('S');

%Sub-Program
function dy=tem14(t,y)
global c1 c2 alpha1 alpha2 eta k m M L g
y1=y(1)-2;
s1=y(2)+c1*y1;s2=y(4)+c2*y(3);
```

```
s=alpha1*s1+alpha2*s2;

f1=(m*L*y(4)^2*sin(y(3))+m*g*sin(y(3))*cos(y(3)))/(M+m*(sin(y(3)))^2);

b1=1/(M+m*(sin(y(3)))^2);

f2=-((m+M)*g*sin(y(3))+m*L*y(4)^2*sin(y(3))*cos(y(3)))/((M+m*(sin(y(3)))^2)*L);

b2=-cos(y(3))/((M+m*(sin(y(3)))^2)*L);

ueq1=-(f1+c1*y(2))/b1;

ueq2=-(f2+c2*y(4))/b2;

usw=-1/(alpha1*b1+alpha2*b2)*(alpha1*b1*ueq2+alpha2*b2*ueq1+eta*sign(s)+k*s);

u=ueq1+ueq2+usw;

dy=[y(2);f1+b1*u;y(4);f2+b2*u];
```